\documentclass[oneside]{amsart}
\usepackage{amsmath,amssymb,color,amsthm,graphicx,cite}
\usepackage{color,bm}
\graphicspath{{../img/}}
\newtheorem{theorem}{Theorem}[section]
\newtheorem{proposition}[theorem]{Proposition}
\newtheorem{lemma}[theorem]{Lemma}
\newtheorem{corollary}[theorem]{Corollary}

\newtheorem{example}{Example}
\theoremstyle{definition}
\newtheorem{definition}{Definition}

\def\Z{\mathbb{Z} }
\def\R{\mathbb{R} }

\def\T{\mathbb{T} }

\def\F{\mathcal{F} }

\author{Tomoo Yokoyama}
\date{\today}
\address{Applied Mathematics and Physics Division, Gifu University, Yanagido 1-1, Gifu, 501-1193, Japan\\}
\email{tomoo@gifu-u.ac.jp}
\thanks{The author was partially supported by JSPS Grant Number 20K03583}

\makeatletter
\@namedef{subjclassname@2020}{%
  \textup{2020} Mathematics Subject Classification}
\makeatother


\title[Morse hyper-graphs and abstract weak element spaces]{Morse hyper-graphs and abstract weak element spaces of semi-decompositions}
\begin{document}
\maketitle

\begin{abstract}
We introduce topological invariants of semi-decompositions (e.g. filtrations, semi-group actions, multi-valued dynamical systems, combinatorial dynamical systems) on a topological space to analyze semi-decompositions from a dynamical systems point of view. In fact, we construct Morse hyper-graphs and abstract weak elements of semi-decompositions. Moreover, the Morse hyper-graphs of the set of sublevel sets of a Morse function of a compact manifold is the Reeb graph of such function as abstract multi-graphs. The abstract weak element space for a simplicial complex is the face poset. The abstract weak element spaces of positive orbits of acyclic directed graphs are their abstract directed multi-graphs. 
\end{abstract}


\section{Introduction}



The concept of a recurrent point is introduced by Birkhoff to describe the limit behavior of orbits \cite{Birkhoff}. 
Conley introduced a weak form of recurrence, called chain recurrence \cite{conley1978isolated}, and showed that dynamical systems on compact metric spaces can be decomposed into blocks, each of which is a chain recurrent one or a gradient one \cite{Conley1988}.
Then this decomposition is called the Morse decomposition, and implies a directed graph, called a Morse graph, which can capture the gradient behaviors. 
The Morse decompositions are developed for mappings and semiflows on arbitrary metric spaces \cite{franks1988pb,hurley1991chain,hurley1992noncompact,hurley1995chain}, set-valued dynamical systems 
\cite{caraballo2015morse,da2016morse,li2007morse,mcgehee1992attractors,wang2020morse},  semigroups and set-valued semiflows \cite{barros2010attractors,barros2012dynamic,bortolan2013skew,barros2010finest,da2017morse,patrao2007morse,patrao2007semiflows,rybakowski2012homotopy}, random dynamical systems \cite{caraballo2012gradient,crauel2004towards,ju2018strong,lin2021morse,liu2005random,liu2007random,liu2007random3,liu2008attractor}, nonautonomous set-valued dynamical systems \cite{aragao2013non,wang2012morse,wang2014morse}, and combinatorial dynamical systems on simplicial complexes \cite{bogdan2020link,mrozek2017conley,tamal2019pers}. 
%
In \cite{yokoyama2021refinements}, the Morse graphs of dynamical systems are refined into abstract orbit spaces. 
The abstract orbit spaces are also refinements of the Reeb graphs of Hamiltonian flows with finitely many singular points on surfaces and the CW decompositions consisting of the unstable manifolds of singular points for Morse flows on closed manifolds. 
Moreover, using the abstract orbit spaces, the author reconstructed some classes of flows from the time-one maps. 
In \cite{yokoyama2021abstract}, the author generalizes such an abstract orbit structure by using concepts of dynamical systems and defines Morse hyper-graphs for decompositions on topological spaces. 
In particular, the abstract element spaces of decompositions are quotient spaces of cell complexes, Morse decompositions, and Reeb graphs. 

Though the studies of decomposition can be applied to foliated spaces and group actions, they can not be applied to semi-group actions in general because the set of orbits of a semi-group need not be a decomposition of the base space. 
To study semi-groups actions, filtrations, and decompositions in a unified manner, a concept of semi-decomposition is introduced in \cite{yokoyama2022relations}.

This paper generalizes abstract (weak) elements and Morse hyper-graphs for decompositions to those for semi-decompositions. 
In particular, we construct topological invariants (e.g. Morse hyper-graphs, abstract elements) of filtrations (e.g. the set of sublevel sets of a function, the set of positive orbits of a multi-valued self-mapping, the set of orbits of a semi-group action (e.g. self-map, semi-flow)) on (possibly non-metrizable) topological spaces. 
Moreover, the abstract element space of a semi-decomposition is a natural generalization of the set of sublevel sets of a Morse function. 

The present paper consists of four sections.
In the next section, recall some concepts and define abstract (weak) element spaces and Morse hyper-graphs for semi-decompositions. 
In \S~3, we demonstrate the existence of Morse hyper-graphs for any semi-decompositions and describe the properties. 
Moreover, we show that the abstract weak element space for a simplicial complex is isomorphic to the face poset, and that the abstract weak element space of positive orbits of an acyclic directed graph is its abstract directed multi-graphs. 
In the final section, some abstract weak element spaces of semi-decompositions are described.

\section{Preliminaries for the description of topological and combinatorial concepts}

In this section, we introduce topological invariants for semi-decompositions of topological spaces. 
To construct such topological invariants, we define some concepts (e.g. Morse graph, recurrence, ``chain recurrence'', ``abstract element'') for semi-decompositions of topological spaces from a dynamical system's point of view. 

\subsection{Fundamental notions of topology and combinatorics}

\subsubsection{Concepts of topology}

For a function $f \colon  X \to \R$ on a topological space $X$, the {\bf Reeb graph} of a function $f \colon  X \to \R$ on $X$ is a quotient space $X/\sim_{\mathrm{Reeb}}$ defined by $x \sim_{\mathrm{Reeb}} y$ if there are a number $c \in \R$ and a connected component of $f^{-1}(c)$ which contains $x$ and $y$.
Notice that the Reeb graph of a Morse function (or, more generally, a function with finitely many critical points) on a closed manifold is a finite graph (see \cite[Theorem~3.1]{saeki2020reeb} for details). 
The inverse image of a value of $\R$ is called the {\bf level set}. 
The inverse image of the half-closed interval $\R_{\leq c}$ for some $c \in \R$ is called the {\bf sublevel set}.

\subsubsection{Concepts of combinatorics}

An {\bf abstract multi-graph} is a triple of sets $V, E$ and a mapping $A \colon E \to \{ \{ x, y \} \mid x,y \in V\}$. 
By a {\bf graph}, we mean a cell complex whose dimension is at most one and which is a geometric realization of an abstract multi-graph.
An {\bf abstract directed multi-graph} is a triple of sets $V, D$ and a mapping $A \colon D \to \{ ( x, y) \mid x,y \in V\}$. 
By a {\bf directed graph}, we mean a cell complex whose dimension is at most one and which is a geometric realization of an abstract directed multi-graph.
A {\bf (abstract) hyper-multi-graph} is a triple of sets $V, E$ and a mapping $A \colon E \to V^*$, where $V^*$ is the family of non-empty finite subsets of $V$. 
An {\bf (abstract) hyper-graph} is a pair of a set $V$ and a family $H \subseteq V^*$.

By a {\bf decomposition}, we mean a family $\mathcal{F}$ of pairwise disjoint nonempty subsets of a set $X$ such that $X = \bigsqcup \mathcal{F}$, where $\bigsqcup$ denotes a disjoint union.
Since connectivity is not required, the sets of orbits of homeomorphisms are also decompositions.


A binary relation $\leq$ on a set $X$ is a {\bf pre-order} if it is reflexive (i.e. $a \leq a$ for any $a \in X$) and transitive (i.e. $a \leq c$ for any $a, b, c \in X$ with $a \leq b$ and $b \leq c$).
For a pre-order $\leq$, the inequality $a<b$ means both $a \leq b$ and $a \neq b$.
A pre-order order $\leq$ is a {\bf total order} (or linear order) if either $a < b$ or $b < a$ for any  points $a \neq b$.
A {\bf chain} is a totally ordered subset of a pre-ordered set with respect to the induced order.
Let $(X, \leq)$ be a pre-ordered set.

\subsubsection{Specialization pre-order for a topological space}

For a topological space $X$, define the {\bf specialization pre-order} $\leq$ by $x \leq y$ if $x \in \overline{\{ y \}}$ for any $x,y \in X$.

\subsection{Abstract (weak) elements for semi-decompositions}

Recall concepts of semi-decompositions and introduce topological invariants for semi-decompositions on topological spaces. 


\subsubsection{Semi-decomposition}

Though the studies of decomposition can be applied to foliated spaces and group actions, they can not be applied to semi-group actions in general because the set of orbits of a semi-group need not be a decomposition of the base space. 
A concept of semi-decomposition is introduced to analyze semi-group actions and decompositions in a unified manner as follows \cite{yokoyama2022relations}.

\begin{definition}
A mapping $\F \colon X \to 2^X$ from a set $X$ is a {\bf semi-decomposition} if it satisfies the following two conditions: 
\\
{\rm(1)} $x \in \F(x)$ for any $x \in X$. 
\\
{\rm(2)} $\F(x) \subseteq \F(y)$ for any $x, y \in X$ with $x \in \F(y)$. 
\end{definition}
By definition, any decompositions are semi-decompositions. 
Note that the set of orbits of any semi-group action is a semi-decomposition and that the set of simplices of the simplicial complex has a semi-decomposition structure. 
Moreover, any filtrations of topological spaces and any stratifications on topological spaces have semi-decomposition structures. 
Here a stratification on a topological space $X$ is a continuous mapping $\pi \colon X \to P$ from $X$ to a poset $P$ with the Alexandroff topology (cf. \cite[Definition A.5.2]{luriehigher}. 

The closure of the element containing a point is called the {\bf element closure} of the point. 
From now on, we use $\F$ for a semi-decomposition on a topological space $X$ unless otherwise stated. 

\subsubsection{A one-to-one correspondence between semi-decompositions and pre-orders}

There is a one-to-one correspondence between semi-decompositions and pre-orders as follows. 
For a semi-decomposition $\F$ on a set $X$, define a pre-order $\leq$ on $X$ by $x \leq y$ if $x \in \F(y)$ (equivalently $\F(x) \subseteq \F(y)$). 
In other words, the induced pre-order is the order induced by the inclusion order on the semi-decomposition $\F$. 
For a pre-order $\leq$ on a set $X$, define a semi-decomposition $\F$ by $\F(y) := \mathop{\downarrow} y = \{ x \in X \mid x \leq y \}$ (i.e. $\F = \{ \mathop{\downarrow} y \mid y \in X \}$). 
In other words, the induced semi-decomposition is the set of downward closures of points.

\subsubsection{Decomposition spaces}
Let $\mathcal{F}$ be a semi-decomposition on a topological space $X$. 
Define the {\bf decomposition space} $X/\F$ as a quotient space $X/\sim_{\F}$, where the equivalence relation $\sim_{\F}$ is defined by $x \sim_{\F} y$ if $\F(x) = \F(y)$.

\subsubsection{Class decompositions}

For any point $x \in X$, define the {\bf element class} $\hat{\mathcal{F}}(x)$ of $x$ as follows:
\[
\hat{\mathcal{F}}(x) := \{ y \in X \mid \overline{\F(y)} = \overline{\F(x)}  \} \subseteq \overline{\F(x)}
\]
Notice that $\F(x) \not\subseteq \hat{\mathcal{F}}(x)$ in general. 
The set $\hat{\mathcal{F}} := \{ \hat{\mathcal{F}}(x) \mid x \in X \}$ is called the {\bf class decomposition} of $\F$. 

\subsubsection{Class semi-decompositions}

To define the element transitive closures for semi-decomposition, we define a binary relation $\leq_{\widetilde{\F},0}$ and a pre-order $\leq_{\widetilde{\F}}$ for a semi-decomposition $\F$ on a topological space $X$ as follows: For any points $x, x'' \in X$, 
\[
x'' \leq_{\widetilde{\F},0} x \text{ if there is a point } x' \in \F(x)  \text{ and } \overline{\mathcal{F}(x'')} = \overline{\mathcal{F}(x')}.
\]
Denote by $\leq_{\F}$ the transitive closure of $\leq_{\widetilde{\F},0}$. 
In other words, we have that $x \leq_{\widetilde{\F}} y$ if there are finitely many points $x_0, x_1, \ldots , x_k \in X$ such that $x_0 = x$, $x_k = y$, and $x_{i-1} \leq_{\widetilde{\F},0} x_{i}$ for any $i \in \{ 1,2, \ldots , k \}$. 
Then the transitive closure $\leq_{\F}$ is a pre-order and contains the specialization pre-order of $\F$ as a subset of $X \times X$. 
For any point $x \in X$, we define the {\bf element transitive closure} $\widetilde{\mathcal{F}}(x)$ of $x$ as follows:
\[
\widetilde{\mathcal{F}}(x) := 
\{ y \in X \mid y \leq_{\widetilde{\F}} x \} \subseteq \overline{\F(x)}
\]
The set $\widetilde{\mathcal{F}} := \{ \widetilde{\mathcal{F}}(x) \mid x \in X \}$ is called the {\bf class semi-decomposition} of $\F$. 
Define an equivalence relation $=_{\widetilde{\F}}$ by $x =_{\widetilde{\F}} y$ if $y \leq_{\widetilde{\F}} x$ and $x \leq_{\widetilde{\F}} y$.  
We have the following observations by the transitivity of the pre-order $\leq_{\widetilde{\F}}$. 

\begin{lemma}\label{lem:rel01}
Any class semi-decomposition is a semi-decomposition. 
\end{lemma}

\begin{lemma}
If a semi-decomposition $\F$ is a decomposition, then its class semi-decomposition corresponds to the class decomposition of the decomposition $\F$. 
\end{lemma}

\begin{lemma}
The following statements hold for any semi-decomposition $\F$ on a topological space $X$ and for any $x,z \in X$: 
\\
{\rm(1)} $\F(x) \subseteq \widetilde{\mathcal{F}}(x) \subseteq \overline{\F(x)}$. 
\\
{\rm(2)} If $x \leq_{\widetilde{\F}} z$, then $\overline{\F(x)} \subseteq \overline{\F(z)}$. 
\\
{\rm(3)} $\hat{\F}(x) = \{ y \in X \mid y =_{\widetilde{\F}} x \} \subseteq \widetilde{\mathcal{F}}(x)$. 
\\
{\rm(4)} If $\F$ is a decomposition, then $\hat{\F}(x) = \widetilde{\mathcal{F}}(x)$. 
\end{lemma}

\begin{proof}
By definition of $\widetilde{\mathcal{F}}(x)$, the assertion {\rm(1)} holds. 
Suppose that $x \leq_{\widetilde{\F}} z$. 
Then there is a sequence $x = x_0 \leq_{\widetilde{\F},0} x_1 \leq_{\widetilde{\F},0} x_2  \leq_{\widetilde{\F},0} \cdots  \leq_{\widetilde{\F},0} x_k = z$. 
By definition, the relation $x_{i} \leq_{\widetilde{\F},0} x_{i+1}$ implies $\overline{\F(x_i)} \subseteq \overline{\F(x_{i+1})}$. 
Therefore $\overline{\F(x)} \subseteq \overline{\F(z)}$ and so the assertion {\rm(2)} holds. 

By the assertion {\rm(2)}, we obtain $\{ y \in X \mid y =_{\widetilde{\F}} x \} \subseteq \{ y \in X \mid \overline{\F(y)} = \overline{\F(x)} \} = \hat{\F}(x)$. 
Fix a point $z \in \hat{\F}(x)$. 
Then $\overline{\F(x)} = \overline{\F(z)}$ and so $z =_{\widetilde{\F}} x$. 
This means that $\hat{\F}(x) = \{ y \in X \mid y =_{\widetilde{\F}} x \} \subseteq \{ y \in X \mid y \leq_{\widetilde{\F}} x \} = \widetilde{\mathcal{F}}(x)$. 

Suppose that $\F$ is a decomposition.
For any points $x'' \leq_{\widetilde{\F},0} x$, since $\F(x') = \F(x)$ for any $x' \in \F(x)$, we have $\overline{\mathcal{F}(x'')} = \overline{\mathcal{F}(x)}$ and so $\hat{\F}(x'') = \hat{\F}(x)$. 
This implies that $y \in \hat{\F}(y) = \hat{\F}(x)$ for any $y \in \widetilde{\mathcal{F}}(x)$. 
Therefore $\hat{\F}(x) = \widetilde{\mathcal{F}}(x)$. 
\end{proof}

\subsubsection{Properties of semi-decompositions}

Let $\F$ be a semi-decomposition of a topological space $X$. 
The subset $\F(A) : = \bigcup_{x \in A} \F(x)$ for any $A \subseteq X$ is called the {\bf saturation} of $A$. 

\begin{definition}
A subset $A$ of $X$ is {\bf {\rm(}$\bm{\F}$-{\rm)}invariant} {\rm(}or {\bf {\rm(}$\bm{\F}$-{\rm)}saturated}{\rm)} if $A = \F(A)$. 
\end{definition}

Notice that a union of elements is not invariant in general. 

\begin{definition}
An invariant subset $A$ of $X$ is {\bf {\rm(}$\bm{\F}$-{\rm)}minimal} if $A = \overline{\F(x)}$ for any point $x \in A$. 
\end{definition}

\subsubsection{Properties of semi-decompositions}

Let $\F$ be a semi-decomposition of a topological space $X$. 

%
%

\begin{definition}
A point $x \in X$ is {\bf closed} with respect to $\F$ if $\overline{\F(x)}= \F(x)$. 
\end{definition}

\begin{definition}
A point $x \in X$ is {\bf proper} with respect to $\F$ if the derived set $\overline{\F(x)} - \F(x)$ is closed. 
\end{definition}

\begin{definition}
A point $x \in X$ is {\bf recurrent} with respect to $\F$ if $\F(x)$ is either closed or non-proper. 
\end{definition}

Denote by $\bm{\mathop{\mathrm{Cl}}(\F)}$ (resp. $\bm{\mathrm{P}(\mathcal{F})}$, $\bm{\mathrm{R}(\mathcal{F})}$, $\bm{\mathcal{R}(\mathcal{F})}$) the set of closed (resp. non-closed proper, non-closed non-proper, recurrent) points. 
Then $X = \mathop{\mathrm{Cl}}(\F) \sqcup \mathrm{P}(\mathcal{F}) \sqcup \mathrm{R}(\mathcal{F}) = \mathrm{P}(\mathcal{F}) \sqcup \mathcal{R}(\mathcal{F})$. 
%
%
%
%
%
From \cite[Lemma~7.1 and Lemma~7.2]{yokoyama2021abstract}, we define invariance of semi-decompositions as follows. 

\begin{definition}
A semi-decomposition $\F$ is {\bf invariant} if it satisfies the following conditions: 
\\
{\rm(1)} The closure of every $\F$-invariant subset is $\F$-invariant. 
\\
{\rm(2)}  The subsets $\mathop{\mathrm{Cl}}(\F)$, $\mathrm{P}(\F)$, and $\mathrm{R}(\F)$ are $\F$-invariant. 
\\
{\rm(3)} For any point $x \in \mathrm{P}(\F)$, the derived set $\overline{\F(x)} - \F(x)$ is $\F$-invariant.  
\end{definition}

Notice that any decompositions satisfy the conditions {\rm(2)} and {\rm(3)} in the previous definition, and so that any invariant decompositions are invariant as semi-decompositions \cite[Lemma~7.1 and Lemma~7.2]{yokoyama2021abstract}. 
Conversely, a decomposition that is invariant with respect to decomposition is invariant with respect to semi-decomposition because of \cite[Lemma~7.2]{yokoyama2021abstract}.

\subsubsection{Properties of invariant semi-decompositions}

We observe the following statements.

\begin{lemma}\label{lem:union_cl}
Let $\F$ be a semi-decomposition $\F$ on a topological space $X$. 
If $\mathop{\mathrm{Cl}}(\F)$ is $\F$-invariant, then $\mathop{\mathrm{Cl}}(\F)$ is $\hat{\F}$-invariant. 
\end{lemma}

\begin{proof}
For any point $x \in \mathop{\mathrm{Cl}}(\F)$, the $\F$-invariance of $\mathop{\mathrm{Cl}}(\F)$ implies that $\hat{\F}(x) = \{ y \in X \mid \overline{\F(y)} = \overline{\F(x)} = \F(x) \} \subseteq \F(x) \subseteq \mathop{\mathrm{Cl}}(\F)$.
This means that $\mathop{\mathrm{Cl}}(\F)$ is $\hat{\F}$-invariant. 
\end{proof}

\begin{lemma}\label{lem:union_p}
Let $\F$ be a semi-decomposition $\F$ on a topological space $X$. 
Suppose that the derived set $\overline{\F(x)} - \F(x)$ for any point $x \in \mathrm{P}(\F)$ is $\F$-invariant. 
Then $\mathrm{P}(\mathcal{F})$ is $\hat{\F}$-invariant and $\hat{\F}(x) \subseteq \widetilde{\F}(x) = \F(x)$ for any $x \in \mathrm{P}(\mathcal{F})$.
\end{lemma}

\begin{proof}
Fix a point $x \in \mathrm{P}(\mathcal{F})$. 
Lemma~\ref{lem:rel01} implies $\hat{\F}(x) \subseteq \widetilde{\F}(x)$. 
We claim that $\hat{\F}(x) \subseteq \widetilde{\F}(x) = \F(x)$. 
Lemma~\ref{lem:rel01} implies $\hat{\F}(x) \subseteq \widetilde{\F}(x)$. 
Assume that there is a point $y \in \widetilde{\F}(x) - \F(x)$. 
Since the derive set $\overline{\F(x)} - \F(x)$ is closed and $\F$-invariant, by $y \in \widetilde{\F}(x) - \F(x) \subseteq \overline{\F(x)} - \F(x)$, we have $\overline{\F(y)} \subseteq \overline{\F(x)} - \F(x)$. 
From $y \in \widetilde{\F}(x)$, there is a point $x' \in \F(x)$ with $\overline{\F(x')} = \overline{\F(y)}$. 
Then $\emptyset \neq \F(x') \subseteq \overline{\F(x')} = \overline{\F(y)} \subseteq \overline{\F(x)} - \F(x)$, which contradicts $\F(x') \subseteq \F(x)$.  
The $\F$-invariance of $\mathrm{P}(\mathcal{F})$ implies that $\mathrm{P}(\mathcal{F})$ is $\hat{\F}$-invariant. 
\end{proof}

\begin{corollary}\label{lem:union_r}
Let $\F$ be a semi-decomposition $\F$ on a topological space $X$. 
Suppose that the derived set $\overline{\F(x)} - \F(x)$ for any point $x \in \mathrm{P}(\F)$ is $\F$-invariant and that $\mathop{\mathrm{Cl}}(\F)$ is $\F$-invariant. 
Then $\mathrm{R}(\mathcal{F})$ is $\hat{\F}$-invariant. 
\end{corollary}

%
%
%
%
%
%

The previous lemmas imply the following invariance. 

\begin{proposition}\label{lem:saturated}
For any invariant semi-decomposition $\F$, the subsets $\mathop{\mathrm{Cl}}(\F)$, $\mathrm{P}(\F)$, and $\mathrm{R}(\F)$ are $\hat{\F}$-invariant.
\end{proposition}

%
%
%

\subsection{Topological concepts for semi-decompositions from Dynamical systems}

We define topological concepts for semi-decompositions from Dynamical systems. 

\subsubsection{Abstract weak element and abstract element for a semi-decomposition}

Let $\mathcal{F}$ be a semi-decomposition on a topological space $X$ and $\bm{p_{X/\F}} \colon X \to X/\F$ the canonical projection. 
Define the following equivalence relation $\sim_{\langle \rangle}$ on $X$: 
$x \sim_{\langle \rangle} y$ if one of the following conditions holds: 
\\
{\rm(1)} $x, y \in \mathrm{Cl}(\mathcal{F})$ and there are points $x' \in \F(x)$ and $y' \in \F(y)$ such that $p_{X/\F}(x')$ and $p_{X/\F}(x')$ are contained in the same connected component of $p_{X/\F}(\mathrm{Cl}(\mathcal{F}))$. 
\\
{\rm(2)} $x, y \in \mathrm{P}(\mathcal{F})$ and there are points $x' \in \F(x)$ and $y' \in \F(y)$ such that $p_{X/\F}(x')$ and $p_{X/\F}(x')$ are contained in the same connected component of $p_{X/\F}(\{ z \in \mathrm{P}(\mathcal{F}) \mid \overline{\F(x)} - \F(x)  = \overline{\F(y)} - \F(y) = \overline{\F(z)} - \F(z) \})$.  
\\
{\rm(3)} $x, y \in \mathrm{R}(\mathcal{F})$ and there are points $x' \in \F(x)$ and $y' \in \F(y)$ such that $p_{X/\F}(x')$ and $p_{X/\F}(x')$ are contained in the same connected component of $p_{X/\F}(\{ z \in \mathrm{R}(\mathcal{F}) \mid \overline{\F(x)} = \overline{\F(y)} = \overline{\F(z)} \})$. 

Similarly, define the following equivalence relation $\sim_{[]}$ on $X$: 
$x \sim_{[]} y$ if one of the following conditions holds: 
\\
{\rm(1)} $x, y \in \mathrm{Cl}(\mathcal{F})$ and there are points $x' \in \F(x)$ and $y' \in \F(y)$ such that $p_{X/\F}(x')$ and $p_{X/\F}(x')$ are contained in the same connected component of $p_{X/\F}(\{ z \in \mathrm{Cl}(\mathcal{F}) \mid \F(x) \cong \F(y) \cong \F(z) \})$, where $ L \cong L'$ means that $L$ and $L'$ are homeomorphic.  
\\
{\rm(2)} $x, y \in \mathrm{P}(\mathcal{F})$ and there are points $x' \in \F(x)$ and $y' \in \F(y)$ such that $p_{X/\F}(x')$ and $p_{X/\F}(x')$ are contained in the same connected component of $p_{X/\F}(\{ z \in \mathrm{P}(\mathcal{F}) \mid \F(x) \cong \F(y) \cong \F(z), \overline{\F(x)} - \F(x) = \overline{\F(y)} - \F(y) = \overline{\F(z)} - \F(z) \})$. 
\\
{\rm(3)} $x, y \in \mathrm{R}(\mathcal{F})$ and there are points $x' \in \F(x)$ and $y' \in \F(y)$ such that $p_{X/\F}(x')$ and $p_{X/\F}(x')$ are contained in the same connected component of $p_{X/\F}(\{ z \in \mathrm{R}(\mathcal{F}) \mid \F(x) \cong \F(y) \cong \F(z), \overline{\F(x)} = \overline{\F(y)} = \overline{\F(z)} \})$. 

For any point $x \in X$, denote by $[x]$ (resp. $\langle x \rangle$) the equivalence class of $x$ with respect to $\sim_{[]}$ (resp. $\sim_{\langle \rangle}$) and called the {\bf abstract weak element}  (resp. {\bf abstract element}) of $x$. 
By definitions, we obtain $[x] \subseteq \langle x \rangle$ for any $x \in X$. 

\subsubsection{Abstract {\rm(}weak{\rm)} element spaces for semi-decompositions}
Let $\mathcal{F}$ be a semi-decomposition on a topological space $X$. 
Define the {\bf abstract weak element space} $X/[\F]$ as a quotient space $X/\sim_{[]}$. 
Similarly, define the {\bf abstract element space} $X/\langle \F \rangle$ as a quotient space $X/\sim_{\langle \rangle}$. 
Since $[x] \subseteq \langle x \rangle$ for any $x \in X$, the abstract element space is a quotient space of the abstract weak element space.

\subsubsection{Correspondence of abstract {\rm(}weak{\rm)} elements for semi-decompositions and decompositions}

Notice that the concepts of abstract weak element and abstract element coincide with those for decompositions under invariance of elements.  
%
The following statement holds.

\begin{lemma}\label{lem:decomp_abstract}
The following statements hold for a decomposition $\F$ on a topological space $X$ and a point $x \in X$: 
\\
{\rm(1)} The abstract element $\langle x \rangle$ is an invariant subset containing $\F(x)$. 
\\
{\rm(2)} If $x \in \mathrm{Cl}(\mathcal{F})$, then the abstract element $\langle x \rangle$ is the inverse image by $p_{X/\F}$ of the connected component of $\mathrm{Cl}(\mathcal{F})/\F \subseteq X/\F$ containing $\F(x)$. 
\\
{\rm(3)} If $x \in \mathrm{P}(\mathcal{F})$, then the abstract element $\langle x \rangle$ is the inverse image by $p_{X/\F}$ of the connected component of $\{ \F(z) \subseteq \mathrm{P}(\mathcal{F}) \mid \overline{\F(x)} - \F(x) = \overline{\F(z)} - \F(z) \} \subseteq X/\F$ containing $\F(x)$. 
\\
{\rm(4)} If $x \in \mathrm{R}(\mathcal{F})$, then the abstract element $\langle x \rangle$ is the inverse image by $p_{X/\F}$ of the connected component of $\{ \F(z) \subseteq \mathrm{R}(\mathcal{F}) \mid \overline{\F(x)} = \overline{\F(z)} \} \subseteq X/\F$ containing $\F(x)$. 
\end{lemma}

\begin{proof}
Let $x, y$ be points in $X$ with $x \sim_{\langle \rangle} y$. 
Since $\F$ is a decomposition, we have $p_{X/\F}(\mathrm{Cl}(\mathcal{F})) = \{ \F(z) \mid z \in \mathrm{Cl}(\mathcal{F}) \} \subseteq X/\F$, $p_{X/\F}(\mathrm{P}(\mathcal{F})) = \{ \F(z) \mid z \in \mathrm{P}(\mathcal{F}) \} \subseteq X/\F$, $p_{X/\F}(\mathrm{R}(\mathcal{F})) = \{ \F(z) \mid z \in \mathrm{R}(\mathcal{F}) \} \subseteq X/\F$, and $p_{X/\F}(x) = \F(x) \in X/\F$. 

Suppose that $x \in \mathrm{Cl}(\mathcal{F})$. 
Let $C$ be the connected component of $p_{X/\F}(\mathrm{Cl}(\mathcal{F})) = \{ \F(z) \mid z \in \mathrm{Cl}(\mathcal{F}) \}$ containing $\F(x)$. 
Then the abstract element $\langle x \rangle$ is the inverse image $p_{X/\F}^{-1}(C)$ and invariant. 

Suppose that $x \in \mathrm{P}(\mathcal{F})$. 
Then $p_{X/\F}(\{ z \in \mathrm{P}(\mathcal{F}) \mid \overline{\F(x)} - \F(x) = \overline{\F(z)} - \F(z) \}) = \{ \F(z) \subseteq \mathrm{P}(\mathcal{F}) \mid \overline{\F(x)} - \F(x) = \overline{\F(z)} - \F(z) \} \subseteq X/\F$. 
Let $C$ be the connected component of $p_{X/\F}(\{ z \in \mathrm{P}(\mathcal{F}) \mid \overline{\F(x)} - \F(x) = \overline{\F(z)} - \F(z) \}) = \{ \F(z) \subseteq \mathrm{P}(\mathcal{F}) \mid \overline{\F(x)} - \F(x) = \overline{\F(z)} - \F(z) \}$ containing $\F(x)$. 
Then the abstract element $\langle x \rangle$ is the inverse image $p_{X/\F}^{-1}(C)$ and invariant. 

Suppose that $x \in \mathrm{R}(\mathcal{F})$. 
Then $p_{X/\F}(\{ z \in \mathrm{R}(\mathcal{F}) \mid \overline{\F(x)} = \overline{\F(z)} \}) = \{  \F(z) \subseteq \mathrm{R}(\mathcal{F}) \mid \overline{\F(x)} = \overline{\F(z)} \} \subseteq X/\F$. 
Let $C$ be the connected component of $p_{X/\F}(\{ z \in \mathrm{R}(\mathcal{F}) \mid \overline{\F(x)} = \overline{\F(z)} \}) = \{ \F(z) \subseteq \mathrm{R}(\mathcal{F}) \mid \overline{\F(x)} = \overline{\F(z)} \}$ containing $\F(x)$. 
Then the abstract element $\langle x \rangle$ is the inverse image $p_{X/\F}^{-1}(C)$ and invariant. 
\end{proof}

We have the following correspondence because of the previous lemma and \cite[Lemma~7.2 and Lemma~7.4]{yokoyama2021abstract}. 

\begin{corollary}\label{cor:decomp_abstract}
If an invariant semi-decomposition on a topological space is a decomposition, then the abstract elements with respect to semi-decomposition correspond to abstract elements with respect to decomposition. 
\end{corollary}

The connectivity of elements implies the following statement. 

\begin{lemma}\label{lem:equiv_abstract}
Let $\F$ be an invariant semi-decomposition of connected elements on a topological space $X$. 
If $\F$ is a decomposition, then the following property holds: 
\[
  \langle x \rangle = \begin{cases}
      \text{the connected component of } \mathrm{Cl}(\F) \text{ containing } \F(x) & \text{if } x \in \mathrm{Cl}(\F) \\
     \text{the connected component of } \\
     \{ x' \in \mathrm{P}(\mathcal{F}) \mid \overline{\F(x)} - \F(x) = \overline{\F(x')} - \F(x') \}) \\
     \text{containing } \F(x) & \text{if } x \in \mathrm{P}(\mathcal{F})\\
     \text{the connected component of} \\
    \{ x' \in \mathrm{R}(\mathcal{F}) \mid \overline{\F(x)} = \overline{\F(x')} \} \text{ containing } \F(x) & \text{if } x \in \mathrm{R}(\mathcal{F}) 
  \end{cases}
\]
\end{lemma}

\begin{proof}
Suppose that $\F$ is a decomposition. 
From \cite[Lemma~7.2]{yokoyama2021abstract}, the decomposition $\F$ is invariant with respect to decomposition and the quotient map $p_{X/\F} \colon X \to X/\F$ is open. 
By \cite[Theorem~6.1.29]{engelking1977general}, the inverse image of any connected subset by $p_{X/\F}$ is connected. 
By Lemma~\ref{lem:decomp_abstract}, any abstract elements are $\F$-invariant and $\F(x) \subseteq \langle x \rangle$ for any $x \in X$. 
Fix a point $x \in X$.
Let $\langle x \rangle'$ be the connected component of $\mathrm{Cl}(\F)$ (resp. $\{ x' \in \mathrm{P}(\mathcal{F}) \mid \overline{\F(x)} - \F(x) = \overline{\F(x')} - \F(x') \})$, $\{ x' \in \mathrm{R}(\mathcal{F}) \mid \overline{\F(x)} = \overline{\F(x')} \}$) containing $\F(x)$ if $x \in \mathrm{Cl}(\F)$ (resp. $\mathrm{P}(\mathcal{F})$, $\mathrm{R}(\mathcal{F})$). 
By \cite[Theorem~6.1.29]{engelking1977general}, the inverse image of any connected subset by $p_{X/\F}$ is connected. 
This implies that the abstract element $\langle x \rangle$ is connected. 

Suppose that $x \in \mathrm{Cl}(\F)$. 
Since $\F(x) \subseteq \mathrm{Cl}(\F)$ is connected, the abstract element $\langle x \rangle$ is the inverse image of the connected component of $\mathrm{Cl}(\F)/\F$ containing $\F(x)$. 
Let $\widetilde{C} \subseteq X/\F$ be the connected component of $p_{X/\F}(\mathrm{Cl}(\F))$ containing $\F(x)$. 
By \cite[Theorem~6.1.29]{engelking1977general}, the inverse image of any connected subset by $p_{X/\F}$ is connected. 
This implies that the inverse image $\langle x \rangle = p_{X/\F}^{-1}(\widetilde{C})$ is connected. 
Then $\langle x \rangle \subseteq \langle x \rangle'$. 
Since the continuous image of any connected subset is connected, we obtain $p_{X/\F}(\langle x \rangle') \subseteq \widetilde{C}$ and so $\langle x \rangle' \subseteq p_{X/\F}^{-1}(\widetilde{C}) = \langle x \rangle$.
Therefore $\langle x \rangle' = \langle x \rangle$. 

Suppose that $x \in \mathrm{P}(\F)$. 
Since $\F(x) \subseteq \mathrm{P}(\F)$ is connected, the abstract element $\langle x \rangle$ is the connected component of $\{ x' \in \mathrm{P}(\mathcal{F}) \mid \overline{\F(x)} - \F(x) =  \overline{\F(x')} - \F(x') \}/\F$ containing $\F(x)$. 
Let $\widetilde{C} \subseteq X/\F$ be the connected component of $p_{X/\F}(\{ x' \in \mathrm{P}(\mathcal{F}) \mid \overline{\F(x)} - \F(x) = \overline{\F(x')} - \F(x') \})$ containing $\F(x)$. 
By \cite[Theorem~6.1.29]{engelking1977general}, the inverse image $\langle x \rangle = p_{X/\F}^{-1}(\widetilde{C})$ is connected and so $\langle x \rangle \subseteq \langle x \rangle'$. 
Since the continuous image of any connected subset is connected, we obtain 
$p_{X/\F}(\langle x \rangle') \subseteq \widetilde{C}$ and so $\langle x \rangle' \subseteq p_{X/\F}^{-1}(\widetilde{C}) = \langle x \rangle$.
Therefore 
$\langle x \rangle' = \langle x \rangle$.

Suppose that $x \in \mathrm{R}(\F)$. 
Since $\F(x) \subseteq \mathrm{R}(\F)$ is connected, the abstract element $\langle x \rangle$ is the connected component of $\{ x' \in \mathrm{R}(\mathcal{F}) \mid \overline{\F(x)} =  \overline{\F(x')} \}$ containing $x$. 
Let $\widetilde{C} \subseteq X/\F$ be the connected component of $p_{X/\F}(\{ x' \in \mathrm{R}(\mathcal{F}) \mid \overline{\F(x)} = \overline{\F(x')} \})$ containing $\F(x)$. 
By \cite[Theorem~6.1.29]{engelking1977general}, the inverse image $\langle x \rangle = p_{X/\F}^{-1}(\widetilde{C})$ is connected and so $\langle x \rangle' \subseteq \langle x \rangle$. 
Since the continuous image of any connected subset is connected, we obtain 
$p_{X/\F}(\langle x \rangle') \subseteq \widetilde{C}$ and so $\langle x \rangle' \subseteq p_{X/\F}^{-1}(\widetilde{C}) = \langle x \rangle$.
Therefore 
$\langle x \rangle' = \langle x \rangle$. 
\end{proof}

The previous lemma implies the following statement.

\begin{lemma}\label{lem:equiv_abstract02}
Let $\F$ be an invariant semi-decomposition of connected elements on a topological space $X$. 
If $\F$ is a decomposition, then the abstract weak elements are abstract elements with respect to decomposition. 
Moreover, the following property holds: 
\[
  [ x ] = \begin{cases}
      \text{the connected component of} \\
   \{ y \in \mathrm{Cl}(\mathcal{F}) \mid \F(x) \cong \F(y) \} \text{ containing } \F(x) & \text{if } x \in \mathrm{Cl}(\F) \\
     \text{the connected component of} \\
      \{ y \in \mathrm{P}(\mathcal{F}) \mid \F(x) \cong \F(y), \overline{\F(x)} - \F(x) = \overline{\F(y)} - \F(y) \} \\
     \text{containing } \F(x) & \text{if } x \in \mathrm{P}(\mathcal{F})\\
     \text{the connected component of} \\
    \{ y \in \mathrm{R}(\mathcal{F}) \mid \F(x) \cong \F(y), \overline{\F(x)} = \overline{\F(y)} \}  \text{ containing } \F(x) & \text{if } x \in \mathrm{R}(\mathcal{F}) 
  \end{cases}
\]
\end{lemma}

\subsection{Morse hyper-graphs of semi-decompositions on topological spaces}

Let $\F$ be a semi-decomposition of a topological space $X$. 

\subsubsection{Quasi-recurrence of points}

\begin{definition}
A point $x \in X$ is {\bf non-maximal} if there is an element $L' \in \F$ such that $\overline{\F(x)} \subsetneq \overline{L'}$. 
\end{definition}

A point is maximal if it is not non-maximal. 
Denote by $\max \F$ the set of maximal points. 
We will define an analogous concept of chain recurrence. 
However, the absence of metrics obstructs a direct interpretation of chain recurrence. 
To construct ``chain recurrence'', recall facts that any points in the orbit closure for a flow are chain recurrent, that any connected components of the chain recurrent point set are unions of abstract elements, and that the chain recurrent set is closed. 
Therefore we define ($\F$-)quasi-recurrence for a topological space, which is the analogous concept of chain recurrence, as follows. 
\begin{definition}
The union of abstract elements of recurrent or non-maximal points is denoted by $\bm{\mathcal{Q}(\F)}$ and called the set of {\bf quasi-recurrent points}. 
\end{definition}

By Lemma~\ref{lem:decomp_abstract}, the union $\mathcal{Q}(\F)$ is $\F$-invariant. 

%

\subsubsection{Morse hyper-graph for a semi-decomposition}

We define the Morse hyper-graph $(V, H)$ for a topological space as follows: 
Let $\F$ be a semi-decomposition on a topological space $X$ with a set $\bm{\mathcal{X}} = \{ \bm{X_{\lambda}} \}_{\lambda \in \Lambda}$ of disjoint subsets $X_\lambda \subseteq X$ ($\lambda \in \Lambda$). 
For any $I \subseteq \lambda$, define a {\bf hyper-edge} $\bm{H_{I}}$ as follows: For any point $x \in X - \bigsqcup_{\lambda \in \Lambda} X_\lambda$, we say that $x \in H_{I}$ if there are disjoint non-empty subsets $(C_i)_{i \in I}$ of the derived set $\overline{\F(x)} - \F(x)$ such that $C_i \subseteq X_i$ and $\overline{\F(x)} - \F(x) = \bigsqcup_{i \in I} C_i$. 
Put $\bm{V} := \{ X_{\lambda} \mid \lambda \in \Lambda \}$ and $\bm{H} := \{ \{ X_i \}_{i \in I} \mid H_{I} \neq \emptyset, I \subseteq \Lambda \}$. 
\begin{definition}
A hyper-graph $\mathcal{G}_{\mathcal{X}} := (V, H)$ is a {\bf Morse hyper-graph} of $\mathcal{X}$ if $X - \bigsqcup_{\lambda \in \Lambda} X_\lambda = \bigsqcup_{I \subseteq \Lambda} H_{I}$. 
\end{definition}

From \cite[Theorem 7.9]{yokoyama2021abstract}, notice that the Morse hyper-graph for an invariant decomposition on a topological space is quotient spaces of the abstract element space and the abstract weak element space with respect to decomposition.
Therefore we define the Morse hyper-graph of a semi-decomposition of a topological space as follows.

\begin{definition}
The hyper-graph $\mathcal{G}_{\mathcal{X}}$ is the {\bf Morse hyper-graph} of $\F$ if $\mathcal{X}$ is the set of inverse images by $q_{\langle  \rangle}$ of connected components of $q_{\langle  \rangle}(\mathcal{Q}(\F))$, where $q_{\langle  \rangle} \colon X \to X/\langle \F \rangle$ is the quotient map. 
\end{definition}

By definition of Morse hyper-graph of a semi-decomposition, the Morse hyper-graph of a topological space is a quotient space of the abstract element space as follows.

\begin{lemma}\label{th:Morse_reduction}
The Morse hyper-graph $\mathcal{G}_{\F}$ for a semi-decomposition $\F$ on a topological space $X$ is a quotient space of the abstract element space $X/\langle \F \rangle$.
\end{lemma}


Lemma~\ref{lem:equiv_abstract} and Corollary~\ref{cor:decomp_abstract} imply the following equivalence. 

\begin{proposition}\label{prop:d_s}
For any invariant decomposition $\F$ on a topological space, the Morse hyper-graphs of $\F$ with respect to decomposition and semi-decomposition coincide with each other. 
 \end{proposition}


The previous proposition implies the following correspondence. 

\begin{corollary}\label{cor:ex_MH_semi_ds}
Let $v \colon \R \times X \to X$ be a flow on a compact metrizable space, $\mathcal{X}$ the set of connected components of the chain recurrent set, and $\F_v$ the set of orbits.  
The Morse hyper-graph $\mathcal{G}_{\mathcal{X}}$ of $\mathcal{X}$ with respect to the {\rm(}semi-{\rm)}decomposition $\F_v$ is the Morse graph $G_{v}$ of $v$ as an undirected graph. 
\end{corollary}


\begin{proof}
Let $\F_v$ be the set of orbits of $v$. 
By Proposition~\ref{prop:d_s}, the Morse hyper-graphs of $\F_v$ with respect to decomposition and semi-decomposition coincide with each other.
Then $\mathop{\mathrm{Cl}}(\F_v) = \mathop{\mathrm{Cl}}(v)$. 
From \cite[Theorem~3.3B]{Conley1988}, the chain recurrent set $\mathop{CR} (v)$ is closed and invariant and contains the non-wandering set $\Omega (v)$, and that connected components of $\mathop{CR} (v)$ are equivalence classes of the relation $\approx_{\mathop{CR}}$ on $\mathop{CR} (v)$, where $x \approx_{\mathop{CR}} y$ if $x \sim_{\mathop{CR}} y$ and $y \sim_{\mathop{CR}} x$.
By \cite[Theorem~6.2]{yokoyama2021abstract}, we have that $\mathcal{R} (\F_v) = \mathcal{R}(v) \subseteq \Omega (v)$. 
Since $\mathcal{R}(\F_v) = \mathop{\mathrm{Cl}}(\F_v) \sqcup \mathrm{R}(\F_v)$ and $\mathcal{R}(v) = \mathop{\mathrm{Cl}}(v) \sqcup \mathrm{R}(v)$, 
the equality $\mathrm{P}(\F_v) = X - \mathcal{R} (\F_v)$ implies that $\mathrm{P}(\F_v) = \mathrm{P}(v)$ and $\mathrm{R}(\F_v) = \mathrm{R}(v)$. 
By definition of non-maximal point, we obtain $\mathcal{R} (\F_v) \cup \max \F_v  \subseteq \Omega (v) \subseteq \mathop{CR} (v)$. 

We claim that $\langle x \rangle \subseteq \mathop{CR} (v)$ for any $x \in \mathop{CR} (v)$. 
Indeed, fix a point $x \in \mathop{CR} (v)$. 
If $x \in \mathcal{R}(v)$, then \cite[Lemma~3.1]{yokoyama2021abstract} implies that $\langle x \rangle \subseteq \mathcal{R} (v) \subseteq \mathop{CR} (v)$. 
Thus we may assume that $x \in \mathrm{P}(v)$. 
Fix any point $y \in \mathrm{P}(v)$ with $\alpha(x) \cup \omega(x) = \alpha(y) \cup \omega(y)$. 
\cite[Lemma~4.1]{yokoyama2021refinements} implies that $\{ x' \in \mathrm{P}(v) \mid \overline{O(x)} - O(x) = \overline{O(x')} - O(x') \neq \emptyset \} = \{ x' \in \mathrm{P}(v) \mid \alpha(x) \cup \omega(x) = \alpha(x') \cup \omega(x') \}$.  
By the closedness of $\mathop{CR} (v)$, since $O(x) \subset \mathop{CR} (v)$, we have $x \sim_{\mathop{CR}} y'$ for any $y' \in \alpha(x) \cup \omega(x) = \alpha(y) \cup \omega(y)$. 
Fix a point $y' \in \alpha(y)$. 
There is a point $y'' \in O(y)$ with $y' \sim_{\mathop{CR}} y''$. 
There is a point $\omega \in \omega(y'') \subseteq \alpha(x) \cup \omega(x)$ with $y'' \sim_{\mathop{CR}} \omega$. 
Since $\omega \in \overline{O(x)} - O(x) \subset \mathop{CR} (v)$, there is a point $x' \in O(x)$ with $y'' \sim_{\mathop{CR}} x'$. 
Since $x,x' \in O(x) \subset \mathop{CR} (v)$, we obtain $x' \sim_{\mathop{CR}} x$ and so $y' \in \mathop{CR} (v)$. 
The invariance of $\mathop{CR} (v)$ implies that $y \in \mathop{CR} (v)$. 
By \cite[Lemma~3.1]{yokoyama2021abstract}, the abstract element $\langle x \rangle$ is a connected component of $\{ x' \in \mathrm{P}(v) \mid \overline{\{ x \}} - \{ x \} = \overline{\{ x' \}} - \{ x' \} \neq \emptyset \} \subseteq \mathop{CR} (v)$. 

Since $\mathcal{R} (\F_v) \cup \max \F_v \subseteq \mathop{CR} (v)$, the previous claim implies that $\mathcal{Q}(\F_v) \subseteq \mathop{CR} (v)$. 
From \cite[Lemma~3.3]{yokoyama2021abstract}, the Morse hyper-graph $\mathcal{G}_\mathcal{X}$ exists. 
By $\bigsqcup_{I \subseteq \Lambda} H_{I} = X - \bigsqcup_{\lambda \in \Lambda} X_\lambda  = X - \mathop{CR} (v) \subseteq \mathrm{P}(v)$, since any $\alpha$-limit sets and $\omega$-limit sets of a flow on a compact Hausdorff space is connected, any hyper-edges $H_I$ are edges, which are non-loops. 
This means that the Morse hyper-graph $\mathcal{G}_{\mathcal{X}}$ is the Morse graph $G_{v}$ as an undirected graph. 
\end{proof}

\section{Existence of Morse hyper-graphs for semi-decompositions and thier properties}

We have the following statements.

\begin{lemma}\label{lem:ch_fol_rec_chain_top}
Let $\F$ be a semi-decomposition on a topological space $X$. 
Then $X - \mathcal{Q}(\F) \subseteq \mathrm{P}(\F) \cap \max \F$. 
\end{lemma}

\begin{proof}
By definition of quasi-recurrence, we have that $X - \mathrm{P}(\F) = \mathcal{R}(\F) \subseteq \mathcal{Q}(\F)$. 
\end{proof}

\begin{lemma}\label{lem:ch_fol_rec_chain_top_02}
The following statements hold for an invariant semi-decomposition $\F$ on a topological space $X$: 
\\
{\rm(1)} $\bigcup_{x \in \mathrm{P}(\F))} (\overline{\F(x)} - \F(x)) \subseteq X - \max \F$. 
\\
{\rm(2)} $\bigcup_{x \in X} \overline{\F(x)} - \F(x) \subseteq \mathcal{Q}(\F)$. 
\\
{\rm(3)} $\mathcal{R}(\F) \sqcup (\bigcup_{x \in X} (\overline{\F(x)} - \F(x)) \cap \mathrm{P}(\F)) \subseteq \mathcal{Q}(\F)$. 
\end{lemma}

\begin{proof}
Fix a point $x \in X$. 
Suppose that $x \in \mathrm{P}(\F)$. 
By invariance of $\F$, the derived set $\overline{\F(x)} - \F(x)$ is $\F$-invariant closed and so 
$\overline{\F(y)} \subseteq \overline{\F(x)} - \F(x) \subsetneq \overline{\F(x)}$ for any $y \in \overline{\F(x)} - \F(x)$. 
Thus the assertion {\rm(1)} holds. 

If $x \in \mathrm{Cl}(\F)$, then $\emptyset = \overline{\F(x)} - \F(x) \subseteq \mathcal{Q}(\F)$. 
Suppose that $x \in \mathrm{R}(\F)$. 
Fix any $y \in \overline{\F(x)} - \F(x)$. 
Then $y \notin \F(x)$. 
If $y \in \mathcal{R}(\F)$, then $y \in \mathcal{Q}(\F)$. 
Thus we may assume that $y \in X - \mathcal{R}(\F) = \mathrm{P}(\F)$. 
Since $\mathcal{P}(\F)$ and $\mathcal{R}(\F)$ are $\F$-invariant, by $x \in \mathcal{R}(\F)$ and $y \in \mathrm{P}(\F)$, we have $\F(x) \cap \F(y) = \emptyset$. 
By invariance of $\F$, the closure $\overline{\F(x)}$ is $\F$-invariant and so $\overline{\F(y)} - \F(y) \subsetneq \overline{\F(y)} \subseteq \overline{\F(x)}$. 
Assume that $x \in \overline{\F(y)}$. 
Then the invariance of $\F$ implies that $\F(x) \subset \overline{\F(y)}$ and so that $\F(x) \subset \overline{\F(y)} - \F(y) \subsetneq \overline{\F(x)}$. 
This means that the closed subset $\overline{\F(y)} - \F(y) \subsetneq \overline{\F(x)}$ contains $\F(x)$, which contradicts the closure $\overline{\F(x)}$ is the minimal closed subset containing $\F(x)$. 
Thus $x \notin \overline{\F(y)}$ and so $\overline{\F(y)} \subsetneq \overline{\F(x)}$. 
This implies that $y$ is non-maximal and so $y \in \mathcal{Q}(\F)$. 

Since $X = \mathcal{R}(\F) \sqcup \mathrm{P}(\F)$, the assertion {\rm(2)} implies the assertion {\rm(3)}. 
%
%
%
%
%
%
%
%
%
%
%
\end{proof}

The following existence of Morse hyper-graph for a topological space holds. 

\begin{proposition}\label{prop:ex_MH_top}
Let be $\F$ a semi-decomposition on a topological space $X$ with the quotient map $p_{X/\F} \colon X \to X/\F$, $\mathcal{Q}' \subseteq X$ a subset containing $\mathcal{Q}(\F)$, and $\mathcal{X} = \{ X_\lambda \}_{\lambda \in \Lambda}$ the set of inverse images by $p_{X/\F}$ of connected components of $\mathcal{Q}'/\F$. 
Then the Morse hyper-graph $\mathcal{G}_{\mathcal{X}}$ exists. 
Moreover, if $\F$ is invariant and $\mathcal{Q}'$ is $\F$-invariant, then the intersection $C_{\lambda} = (\overline{\F(x)} - \F(x)) \cap X_{\lambda}$ are $\F$-invariant. 
\end{proposition}

\begin{proof}
Lemma~\ref{lem:ch_fol_rec_chain_top} implies that $X - \bigcup_{\lambda \in \Lambda} X_\lambda = X - \mathcal{Q}' \subseteq X - \mathcal{Q}(\F) \subseteq  \mathrm{P}(\F) \cap \max \F$. 
Moreover, fix any point $ x \in X - \bigcup_{\lambda \in \Lambda} X_\lambda \subset \mathrm{P}(\F) \cap \max \F$. 
Since $\overline{\F(x)} - \F(x) \subsetneq \overline{\F(x)}$, by Lemma~\ref{lem:ch_fol_rec_chain_top_02}, we have that $\overline{\F(x)} - \F(x) \subseteq X - \max \F \subseteq \mathcal{Q}(\F)$. 
Therefore define an index set $I \subseteq \Lambda$ by $i \in I$ if $C_i := (\overline{\F(x)} - \F(x)) \cap X_i \neq \emptyset$. 
Then $\overline{\F(x)} - \F(x) = \bigsqcup_{i \in I} C_i$ and so $x \in  H_I$. 
This means that $X - \bigcup_{\lambda \in \Lambda} X_\lambda = X - \mathcal{Q}' = \bigsqcup_{I \subseteq \Lambda} H_{I}$. 

Suppose that $\F$ is invariant and $\mathcal{Q}'$ is $\F$-invariant. 
Then any inverse images $X_i$ of $p_{X/\F}$ and the derived set $\overline{\F(x)} - \F(x)$ are $\F$-invariant. 
Therefore the intersections $C_i = (\overline{\F(x)} - \F(x)) \cap X_i$ are $\F$-invariant. 
%
\end{proof}

Since the union $\mathcal{Q}(\F)$ is $\F$-invariant, we have the following statement. 

\begin{corollary}\label{lem:ex_MH_top}
The Morse hyper-graph for any semi-decomposition of a topological space exists.
\end{corollary}

\subsection{Generalization of Reeb graphs of Morse functions}

Recall that the Reeb graph of a Morse function on a closed manifold is the abstract weak element space of the set of connected components of level sets as abstract multi-graphs \cite[Proposition 7.6]{yokoyama2021abstract}. 
For a semi-decomposition of sublevel sets, we have a similar statement. 
In other words, we show that the abstract element space of a semi-decomposition is a natural generalization of the Reeb graph of a Morse function. 

\begin{proposition}\label{prop:reeb}
The Reeb graph of a Morse function on a closed manifold is the abstract weak element space of the set of connected components of sublevel sets as abstract multi-graphs. 
\end{proposition}

\begin{proof}
Let $f$ be a $C^1$ function with finitely many critical points on a closed manifold $M$ and $\F$ the set of connected components of sublevel sets of $f$. 
Then any elements of $\F$ are closed. 
Notice that the Reeb graph of $f$ is a finite graph as abstract graph (cf. \cite[Proposition 7.6]{yokoyama2021abstract}). 
Denote by $\mathrm{Crit}(f)$ the set of critical points.
By Morse theory, for any points $x,y$ contained in the same connected component of the complement $X - \bigcup_{x \in \mathrm{Crit}(f)} \F(x)$, the elements $\F(x)$ and $\F(y)$ are homeomorphic to each other and are submanifolds. 
Therefore Lemma~\ref{lem:equiv_abstract02} implies that any abstract weak elements are connected component of the complement $X - \bigcup_{x \in \mathrm{Crit}(f)} \F(x)$ and connected components of $\bigcup_{x \in \mathrm{Crit}(f)} \F(x)$.  
Since any connected component of the complement $X - \bigcup_{x \in \mathrm{Crit}(f)} \F(x)$ correspond to abstract edges and any connected components of $\bigcup_{x \in \mathrm{Crit}(f)} \F(x)$ correspond to vertices,  the Reeb graph of $f$ is the abstract weak element space $M/\F$ as an abstract multi-graph. 
\end{proof}

\subsection{Generalization of face posets of simplicial complices}

Let $\mathcal{K}$ be a simplicial complex and $|\mathcal{K}| := \bigcup \mathcal{K}$ the polyhedron. 
For any point $x \in |\mathcal{K}|$, denote by $\F(x)$ the lowest dimensional simplex of $\mathcal{K}$ containing $x$. 
Then the set $\bm{\F_{\mathcal{K}}} := \{ \F(x) \mid x \in |\mathcal{K}| \}$ is a semi-decomposition of the polyhedron $|\mathcal{K}|$ and is called the {\bf semi-decomposition associated to the simplicial complex} $\mathcal{K}$. 
The {\bf face poset} of a simplicial complex $\mathcal{K}$ is defined by $(\mathcal{K}, \subseteq)$, where $\subseteq$ is the inclusion order on $|\mathcal{K}|$.
We have the following statement. 

\begin{proposition}
The abstract weak element space $|\mathcal{K}|/\langle \F_{\mathcal{K}} \rangle$ for a simplicial complex $\mathcal{K}$ with the specialization order is isomorphic to the face poset of $\mathcal{K}$. 
Moreover, we have that $|\mathcal{K}|/\langle \F_{\mathcal{K}} \rangle = \bigcup_{k \geq 0} \{ \sigma \setminus |\mathcal{K}^{k-1}| \mid \sigma \in \mathcal{K}^k \}$ as a set, where $\mathcal{K}^i$ is the set of $i$-cells of $\mathcal{K}$. 
\end{proposition}

\begin{proof}
The union $|\mathcal{K}^0|$ of $0$-cells of $\mathcal{K}$ is a set of discrete points and so $\langle x \rangle = \sigma$ for any $0$-cell $\sigma \in \mathcal{K}^0$ and any point $x \in \sigma$. 
For any $1$-cell $\sigma \in \mathcal{K}^1$ and any points $x,y \in \sigma \setminus |\mathcal{K}^0|$, we have that $\F(x) = \F(y) = \sigma$ and so $\langle x \rangle = \langle y \rangle = \sigma \setminus |\mathcal{K}^0|$. 
This means that the abstract weak element of any point $x$ of a $1$-cell $\sigma \in \mathcal{K}^1$ is the interior $\sigma \setminus |\mathcal{K}^0|$ in $|\mathcal{K}^1|$. 
Similarly, for any $k$-cell $\sigma \in \mathcal{K}^k$ and any points $x,y \in \sigma \setminus |\mathcal{K}^{k-1}|$, we have that $\F(x) = \F(y) = \sigma$ and so $\langle x \rangle = \langle y \rangle = \sigma \setminus |\mathcal{K}^{k-1}|$. 
This means that the abstract weak element of any point $x$ of a $k$-cell $\sigma \in \mathcal{K}^k$ is the interior $\sigma \setminus |\mathcal{K}^{k-1}|$ in $|\mathcal{K}^k|$. 
The quotient space $|\mathcal{K}|/\langle \F_{\mathcal{K}} \rangle = \bigcup_{k \geq 0} \{ \sigma \setminus |\mathcal{K}^{k-1}| \mid \sigma \in \mathcal{K}^k \}$ is the set of interiors of cells of $\mathcal{K}$ as a set. 
By construction, for any simplices $\sigma, \mu \in \mathcal{K}$, the simplex $\sigma$ is a face of $\mu$ if and only if $\sigma \subseteq \mu$.
Let $\leq$ be the specialization pre-order on $|\mathcal{K}|/\langle \F_{\mathcal{K}} \rangle$. 
Fix any points $x, y \in |\mathcal{K}|$. 
By definition, we have that $\langle x \rangle \leq \langle y \rangle$ if and only if $\langle x \rangle \in \overline{\{ \langle y \rangle\}}^{|\mathcal{K}|/\langle \F_{\mathcal{K}} \rangle}$ (i.e. the abstract weak element $\langle x \rangle$ is contained in the closure of the singleton of the abstract weak element $\langle y \rangle$ in the abstract weak element space $|\mathcal{K}|/\langle \F_{\mathcal{K}} \rangle$). 
Since the abstract weak element $\langle x \rangle$ (resp. $\langle y \rangle$) is the interior of the cell $\F(x)$ (resp. $\F(y)$) and  the cell $\F(x)$ (resp. $\F(y)$) is the closure $\overline{\langle x \rangle}$ (resp. $\overline{\langle y \rangle}$) of the interior of $\F(x)$ (resp. $\F(y)$), we obtain that $\langle x \rangle \leq \langle y \rangle$ if and only if $\F(x) \subseteq \F(y)$. 
This means that the abstract weak element space $|\mathcal{K}|/\langle \F_{\mathcal{K}} \rangle$ with the specialization order is isomorphic to the face poset of $\mathcal{K}$. 
\end{proof}

\subsection{Generalization of abstract directed multi-graph structures of acyclic directed graphs}

Let $X$ be a directed graph with the vertex set $V$. 
By definition, the directed graph is a cell complex whose dimension is at most one. 
For any point $x \in X - V$, denote by $\bm{d(x)}$ the directed edge containing $x$ and by $\bm{a(x)} \subsetneq d(x)$ the arc from $x$ to the tail of $d(x)$. 
For any $x \in V$, denote by $\bm{X_{\leq 1}(x)}$ the union of $x$, adjacent vertices of $x$, and the directed edges from $x$, and define $\bm{X_{\leq k}(x)} := \bigcup_{y \in V \cap X_{\leq k-1}(x)} X_{\leq 1}(y)$. 
For any $x \in V$, the positive orbit $\bm{O^+(x)}$ is defined by $\bigcup_{k =1}^\infty X_{\leq k}(x)$. 
For any point $x \in X - V$, the positive orbit $\bm{O^+(x)}$ is defined by the union $a(x) \cup O^+(y)$, where $y$ is the tail of $d(x)$. 
This means that the positive orbit of a point $x$ is the union of any arc $C \colon [0,1] \to X$ from $x$ whose orientation on $C([0,1]) \setminus V$ corresponds to the direction of the edges of $X$. 
Then the set of positive orbits is a semi-decomposition, denoted by $\bm{\F_{X}}$ and called the semi-decomposition of positive orbits of the directed graph $X$. 
A directed graph is acyclic if it has no directed cycles. 
Note that any loops are directed cycles. 
We have the following statement. 

\begin{proposition}
The abstract weak element space $X/\langle \F_X \rangle$ of any finite directed graph $X$ is an abstract directed multi-graphs whose directed graph is obtained from $X$ by collapsing any directed cycles. 
Moreover, if $X$ is acyclic, then the abstract weak element space $X/\langle \F_X \rangle$ corresponds to the abstract directed multi-graph of $X$. 
\end{proposition}

\begin{proof}
Let $X$ be a directed graph with the vertex set $V$. 
We may assume that $X$ is non-empty and connected. 
For any directed edge $d$ and any points $x,y \in d \setminus V$, the positive orbits $O^+(x)$ and $O^+(y)$ are homeomorphic to each other. 
Let $\F_X$ be the set of positive orbits.  
For any point $x \in X - V$, the abstract weak element $\langle x \rangle$ of $x$ contains $d(x)$. 
Fix a vertex $x \in V$. 

We claim that we may assume that $X$ has at least two vertices. 
Indeed, assume that $X$ has exactly one vertex. 
Then any positive orbits are $X$ and so $X/\langle \F_X \rangle$ is a singleton. 

We claim that $\langle x \rangle \neq \{ x \}$ if and only if $x$ either is contained in any directed cycles or has exactly one out-degree. 
Indeed, notice that the abstract weak element of a point contained in a directed cycle is not a singleton. 
Suppose that $\langle x \rangle \neq \{ x \}$. 
We may assume that $x$ is not contained in any directed cycles. 
Since $\langle x \rangle \neq \{ x \}$, there is a directed edge $d$ from $x$. 
Assume that $x$ has at least two out-degree. 
For any point $y \in X_{\leq 1}(x) \setminus V$, the finiteness of $X$ implies that the base spaces of the positive orbits $O^+(y) \subsetneq O^+(x)$ are finite graphs and not homeomorphic to each other. 
This means that $\langle x \rangle = \{ x \}$, which contradicts that $\langle x \rangle \neq \{ x \}$. 
Conversely, if $x$ is contained in any directed cycles, then $\langle x \rangle \neq \{ x \}$. 
Suppose that $x$ has exactly one out-degree. 
Then let $d$ be the directed edge from $x$. 
For any point $y \in d \setminus V$, the positive orbits $O^+(y)$ and $O^+(x)$ are homeomorphic to each other. 
This implies that $y \in \langle x \rangle \neq \{ x \}$. 

By the previous claim, the abstract weak element space $X/\langle \F_X \rangle$ is homeomorphic to the abstract weak element space of positive orbits of the resulting finite directed graph from $X$ collapsing directed cycles into singletons. 
Therefore we may assume that $X$ is acyclic by collapsing directed cycles into singletons if necessary. 
Then $\langle x \rangle \neq \{ x \}$ if and only if $x$ has exactly one out-degree. 
This means that the abstract weak element space $X/\langle \F_X \rangle$ is homeomorphic to the abstract weak element space of positive orbits of the resulting finite directed graph from $X$ collapsing edges whose heads have exactly out-degree one into singletons. 
Therefore we may assume that $X$ has no vertices whose out-degrees are one, by collapsing edges whose heads have exactly out-degree one into singletons if necessary. 
Then the abstract weak elements of any vertices are singletons, and so the abstract weak elements of any point $y \in X - V$ are directed edges $d(y)$.  
\end{proof}

\section{Examples of abstract {\rm(}weak{\rm)} element spaces}

We state some examples of abstract (weak) element spaces. 

\begin{example}\label{ex:01}
{\rm
Let $v_0 \colon \R \times \T^2 \to \T^2$ be an irrational rotational flow on a torus $\T^2 := \R^2/\Z^2$ generated by a vector field $V_0([x,y]) = (1, \alpha)$ for some $\alpha \in \R - \mathbb{Q}$. 
The set $\F_0$ of positive orbits $O^+(p) = \{ v_0(t,p) \mid t \in \R_{\geq 0}\}$ is a semi-decomposition. 
Then the abstract (weak) element space $\T^2/[\F_0] = \T^2/\langle \F_0 \rangle$ is a singleton. 
}
\end{example}

\begin{example}
{\rm
Let $V$ be the resulting vector field from $V_0$ by replacing the point $0 := [0,0] \in \T^2$ to a singular point using a bump function. 
Then the resulting flow $v$ generated by $V$ consists of one singular point $0$ and non-closed recurrent orbits. 
Denote by $O_+$ (resp. $O_-$) the orbit of $v$ whose $\alpha$-limit (resp. $\omega$-limit) set is $0$. 
The set $\F$ of positive orbits of $v$ is a semi-decomposition. 
Then $\overline{\F(0)} = \F(0) = \{ 0 \}$. 
For any $p \in \T^2$, we have that  $\overline{\F(p)} = \F(p) \sqcup  \{ 0 \}$ if and only if $p \in O_+$. 
Therefore $\mathrm{Cl}(\mathcal{F}) = \{ 0 \}$, $\mathrm{P}(\mathcal{F}) = O_+$, and $\mathrm{R}(\mathcal{F}) = \T^2 - (O_- \sqcup \{ 0 \}) = \T^2 - \overline{O_-}$. 
This implies that the abstract (weak) element space $\T^2/[\F] = \T^2/\langle \F \rangle$ is a topological space which consists of three elements $\mathrm{Cl}(\mathcal{F})$, $\mathrm{P}(\mathcal{F})$, and $\mathrm{R}(\mathcal{F})$. 
Moreover, the specialization pre-order of $\T^2/[\F]$ is a total order such that $\mathrm{Cl}(\mathcal{F}) < \mathrm{P}(\mathcal{F}) <  \mathrm{R}(\mathcal{F})$ (i.e. $\{ 0 \} < O_+ < \T^2 - \overline{O_-}$). 
}
\end{example}

\bibliographystyle{abbrv}
\bibliography{yt20211011}

\end{document}